\documentclass{amsart}
\usepackage[colorlinks]{hyperref}
\AtBeginDocument{
  \hypersetup{
    linkcolor=blue,
    citecolor=green,
  }
}
\usepackage[displaymath,mathlines,pagewise]{lineno}
\newcommand*\patchAmsMathEnvironmentForLineno[1]{%
  \expandafter\let\csname old#1\expandafter\endcsname\csname #1\endcsname
  \expandafter\let\csname oldend#1\expandafter\endcsname\csname end#1\endcsname
  \renewenvironment{#1}%
     {\linenomath\csname old#1\endcsname}%
     {\csname oldend#1\endcsname\endlinenomath}}%
\newcommand*\patchBothAmsMathEnvironmentsForLineno[1]{%
  \patchAmsMathEnvironmentForLineno{#1}%
  \patchAmsMathEnvironmentForLineno{#1*}}%
\AtBeginDocument{%
\patchBothAmsMathEnvironmentsForLineno{equation}%
\patchBothAmsMathEnvironmentsForLineno{align}%
\patchBothAmsMathEnvironmentsForLineno{flalign}%
\patchBothAmsMathEnvironmentsForLineno{alignat}%
\patchBothAmsMathEnvironmentsForLineno{gather}%
\patchBothAmsMathEnvironmentsForLineno{multline}%
}

\usepackage{cleveref}
\crefname{equation}{}{}

\newtheorem{theorem}{Theorem}[section]
\newtheorem{lemma}[theorem]{Lemma}
\newtheorem{corollary}[theorem]{Corollary}

\theoremstyle{definition}

\theoremstyle{remark}
\newtheorem{remark}[theorem]{Remark}

\numberwithin{equation}{section}

\begin{document}

\title{A note on the BMO and Calder\'{o}n-Zygmund estimate}

\author{Yuanyuan Lian}
\address{School of Mathematical Sciences, Shanghai Jiao Tong University, Shanghai, China}
\email{lianyuanyuan@sjtu.edu.cn; lianyuanyuan.hthk@gmail.com}

\author{Kai Zhang}
\address{School of Mathematical Sciences, Shanghai Jiao Tong University, Shanghai, China}
\email{zhangkaizfz@gmail.com}
\thanks{This research is supported by the China Postdoctoral Science Foundation (Grant No.
2021M692086), the National Natural Science Foundation of China (Grant No. 12031012 and 11831003) and  the Institute of Modern Analysis-A Frontier Research Center of Shanghai.}


\date{March 24, 2021}



\begin{abstract}
In this note, we give a simple proof of the pointwise BMO estimate for Poisson's equation. Then the Calder\'{o}n-Zygmund estimate follows by the interpolation and duality.
\end{abstract}

\maketitle

\section{Introduction}\label{S1}
Consider the following Poisson's equation
\begin{equation}\label{e1.1}
\Delta u= f ~~~~\mbox{in}~B_1,
\end{equation}
where $B_1\subset R^n$ is the unit ball. In 1952, Calder\'{o}n and Zygmund (\cite{MR0052553}, see also \cite[Chapter 9]{MR1814364}) proved the classical $W^{2, p}$ estimate ($1<p<\infty$) for \cref{e1.1} by the method of singular integral. That is, for any strong solution $u\in W^{2, p}(B_1)$ of \cref{e1.1},
\begin{equation*}
  \|u\|_{W^{2,p}(B_{1/2})}\leq C(\|u\|_{L^p(B_1)}+\|f\|_{L^{p}(B_1)}),
\end{equation*}
where $C$ depends only on $n$ and $p$. In 2003, Wang \cite{MR1987802} gave a direct and elementary proof of the $W^{2, p}$ estimate without using the theory of singular integral. Wu, Yin and Wang \cite[Chapter 9]{MR2309679} presented another elementary proof based on energy estimates.

It is also well known that the $W^{2, p}$ estimate fails for $p=1$ and $p=\infty$ (see \cite[Chapter 7.1.3]{MR3099262}). It was found that the Hardy space $H^1$ and the BMO (bounded mean oscillation) space are appropriate substitutes for $L^1$ and $L^{\infty}$ respectively. The estimates in the $H^1$ and the BMO space are proved by means of singular integral as well (see \cite[Theorem 4, Chapter 3]{MR1232192} and \cite{MR209917}).

In this note, we give a simple proof of the $W^{2, BMO}$ estimate for \cref{e1.1} (see \Cref{t-2} below). By combining with the interpolation and duality argument, we have the $W^{2,p}$ estimate for any $1<p<\infty$.

We introduce some notations. The $B_r(x)\subset R^n$ denotes the ball with radius $r$ and center $x$, and $B_r:=B_r(0)$. Let
\begin{equation*}
f_{\Omega}:= \frac{1}{|\Omega|}\int_{\Omega} f~~\mbox{and}~~\|f\|^*_{L^p(\Omega)}:= \left(\frac{1}{|\Omega|}\int_{\Omega} |f|^p\right)^{1/p}, ~\forall ~1\leq p<\infty,
\end{equation*}
where $\Omega\subset R^n$ is a bounded domain and $|\Omega|$ denotes its Lebesgue measure.

For $x_0\in \Omega$ and $r_0>0$, if
\begin{equation*}
   |f|_{*,x_0}:=\sup_{0<r<r_0} \| f-f_{B_r(x_0)\cap \Omega}\|^*_{L^2(B_r(x_0)\cap \Omega)}<+\infty,
\end{equation*}
we say that $f$ is BMO at $x_0$ or $f\in BMO(x_0)$ (with radius $r_0$). If $f\in BMO(x_0)$ for any $x_0\in \Omega$ with the same radius $r_0$ and
\begin{equation*}
 |f|_{*,\Omega}:=\sup_{x\in \Omega}|f|_{*,x} <+\infty,
\end{equation*}
we say that $f$ is a BMO function or $f\in BMO(\Omega)$ (with radius $r_0$). Moreover, we endow $BMO(\Omega)$ with the following norm:
\begin{equation*}
  \|f\|_{BMO(\Omega)}:=\|f\|_{L^2(\Omega)}+|f|_{*,\Omega}.
\end{equation*}
If $f$ has weak derivatives (denoted by $Df$) and $Df\in BMO(x_0)$ ($BMO(\Omega)$), we say that $u\in W^{1, BMO}(x_0)$ ($W^{1, BMO}(\Omega)$). Similarly, we can define $W^{2, BMO}$ etc.

Note that the usual BMO space is defined based on $L^1$ norm rather than $L^2$ norm. In fact, these definitions are equivalent (see \cite[Corollary on P. 144 ]{MR1232192} and \cite[Corollary 6.22]{MR3099262}).

Unless stated otherwise, $C$ always denotes a constant depending only on the dimension $n$ throughout this note.

\section{Main result}
Before starting to prove the main result, we describe the idea briefly. We adopt perturbation argument, i.e., we use harmonic functions to approximate the solution of \cref{e1.1}. If $f$ is small enough, then $u$ can be approximated by a polynomial of degree $2$ at some scale (called key step, see \Cref{l-1}). This is proved by the method of compactness. Next, by the standard scaling argument, we have a sequence of polynomials approximating $u$ at different scales (see \Cref{t-1}). Then this essentially implies that $u\in BMO(0)$ (see \Cref{t-2}). This technique has been used widely to prove the $C^{k,\alpha}$ regularity since the seminal work of Caffarelli (see \cite{MR1005611, MR1351007, zbMATH07405906,MR4088470, MR3246039, MR1139064} etc.). It turns out that the BMO regularity can be regarded as a kind of pointwise regularity as the $C^{k,\alpha}$ regularity.

First, we prove the key step.
\begin{lemma}\label{l-1}
Let $u\in W_{loc}^{1,2}(B_1)$ be a weak solution of
\begin{equation}\label{e1}
\Delta u=f ~~~~\mbox{in}~B_1.
\end{equation}
Suppose that $\|u\|^*_{L^2(B_1)}\leq 1$ and $\|f\|^*_{L^2(B_1)}\leq \delta$, where $0<\delta<1$ depends only on $n$.

Then there exists a polynomial $P$ of degree $2$ such that
\begin{equation}
  \begin{aligned}
    &\|u-P\|^*_{L^2(B_{\eta})}\leq \eta^2,\\
    &\Delta P=f_{B_{\eta}},\\
    &|P(0)|+|DP(0)|+|D^2P(0)|\leq \bar{C},
  \end{aligned}
\end{equation}
where $0<\eta<1$ and $\bar{C}$ depend only on $n$.
\end{lemma}
\proof We prove the lemma by contradiction. Suppose that the lemma is false. Then there exist sequences of $u_m$ and $f_m$ such that
\begin{equation*}
  \Delta u_m=f_m ~~~~\mbox{in}~B_1
\end{equation*}
with $\|u_m\|^*_{L^2(B_1)}\leq 1$ and $\|f_m\|^*_{L^2(B_1)}\leq 1/m$. But for any polynomial $P$ of degree $2$ with $\Delta P=f_{m, B_{\eta}}$ and $|P(0)|+|DP(0)|+|D^2P(0)|\leq \bar{C}$, we have
\begin{equation}\label{e2}
\|u_m-P\|^*_{L^2(B_{\eta})}> \eta^2,
\end{equation}
where $0<\eta<1$ and $\bar{C}$ are to be specified later.

Then there exist $u\in W_{loc}^{1,2}(B_1)$ and subsequences of $u_m$ (denoted by $u_m$ again) such that $u_m\rightarrow u$ in $W_{loc}^{1,2}(B_1)$ weakly and in $L^2_{loc} (B_1)$ strongly. Moreover,
\begin{equation*}
  \Delta u=0~~~~\mbox{in}~B_1.
\end{equation*}
Since $u$ is a harmonic function, there exists a polynomial $\bar{P}$ of degree $2$ such that
\begin{equation*}
  \begin{aligned}
    &\|u-\bar{P}\|^*_{L^2(B_{r})}\leq Cr^3, ~\forall ~0<r<1/2,\\
    &\Delta \bar{P}=0,\\
    &|\bar{P}(0)|+|D\bar{P}(0)|+|D^2\bar{P}(0)|\leq C.
  \end{aligned}
\end{equation*}
Take $\bar{C}=C+1$ and $\eta$ small enough such that
\begin{equation*}
  \eta \bar{C}\leq 1/2.
\end{equation*}
Thus,
\begin{equation}\label{e3}
 \|u-\bar{P}\|^*_{L^2(B_{\eta})}\leq \frac{1}{2}\eta^2.
\end{equation}

Set
\begin{equation*}
P_m(x)=\bar{P}(x)+\frac{f_{m, B_{\eta}}}{2n}|x|^2.
\end{equation*}
Then $\Delta P_m=f_{m, B_{\eta}}$ and $|P_m(0)|+|DP_m(0)|+|D^2P_m(0)|\leq \bar{C}$ (for $m$ large enough). By \cref{e2},
\begin{equation*}
\|u_m-P_m\|^*_{L^2(B_{\eta})}> \eta^2.
\end{equation*}
Let $m\rightarrow \infty$ (noting $f_{m, B_{\eta}}\rightarrow 0$) and we have
\begin{equation*}
\|u-\bar{P}\|^*_{L^2(B_{\eta})}\geq \eta^2,
\end{equation*}
which contradicts with \cref{e3}. \qed~\\

Now, we give the scaling argument.
\begin{lemma}\label{t-1}
Let $0<\delta<1$ be as in \Cref{l-1} and $u\in W^{1,2}_{loc}(B_1)$ be a weak solution of
\begin{equation*}
\Delta u=f ~~~~\mbox{in}~B_1.
\end{equation*}
Suppose that $\|u\|^*_{L^2(B_1)}\leq 1$, $\|f\|^*_{L^2(B_1)}\leq \delta$ and $|f|_{*,0}\leq \delta$ (i.e., $\|f-f_{B_r}\|^*_{L^2(B_r)}\leq \delta$ for any $0<r<1$).

Then there exist a sequence of polynomials $P_m$ of degree $2$ such that for any $m\geq 1$,
\begin{equation}\label{e4}
  \begin{aligned}
    &\|u-P_m\|^*_{L^2(B_{\eta^m})}\leq \eta^{2m},\\
    &\Delta P_m=f_{B_{\eta^m}},\\
    &|(P_m-P_{m-1})(0)|+\eta^{m-1}|D(P_m-P_{m-1})(0)|+\eta^{2(m-1)}|D^2(P_m-P_{m-1})(0)|\leq \bar{C}\eta^{2(m-1)},
  \end{aligned}
\end{equation}
where $0<\eta<1$ and $\bar{C}$ are as in \Cref{l-1}.
\end{lemma}
\proof We prove the lemma by induction. For $m=1$, by setting $P_0\equiv 0$ and \Cref{l-1}, the conclusion holds. Suppose that the conclusion holds for $m$ and we need to prove that it holds for $m+1$.

Let $r=\eta^m$, $y=x/r$ and
\begin{equation*}
  v(y)=\frac{u(x)-P_m(x)}{r^2}.
\end{equation*}
Then
\begin{equation*}
  \Delta v=\tilde{f}~~~~\mbox{in}~B_1,
\end{equation*}
where $\tilde{f}(y)=f(x)-\Delta P_m=f(x)-f_{B_r}$. Clearly, $\|v\|^*_{L^2(B_1)}\leq 1$. In addition,
\begin{equation*}
  \|\tilde{f}\|^*_{L^2(B_1)}=\|f-f_{B_r}\|^*_{L^2(B_r)}\leq \delta.
\end{equation*}
By \Cref{l-1}, there exists a polynomial $P$ of degree $2$ such that
\begin{equation*}
  \begin{aligned}
    &\|v-P\|^*_{L^2(B_{\eta})}\leq \eta^2,\\
    &\Delta P=\tilde{f}_{B_{\eta}},\\
    &|P(0)|+|DP(0)|+|D^2P(0)|\leq \bar{C}.
  \end{aligned}
\end{equation*}
By rescaling back to $u$ with $P_{m+1}(x)=P_m(x)+\eta^{2k}P(y)$, \cref{e4} holds for $m+1$. By induction, the proof is complete.\qed~\\

Now, we show that \Cref{t-1} implies the $W^{2, BMO}$ regularity of the solution.
\begin{theorem}\label{t-2}
Let $u\in W^{1,2}_{loc}(B_1)$ be a weak solution of
\begin{equation*}
\Delta u=f ~~~~\mbox{in}~B_1.
\end{equation*}
Suppose that $f\in BMO(0)$ with radius $1$. Then $u\in W^{2, BMO}(0)$ with radius $\eta^2$ and
\begin{equation*}
  |D^2u|_{*,0}\leq C\left(\|u\|_{L^2(B_1)}+\|f\|_{L^2(B_1)}+|f|_{*,0}\right),
\end{equation*}
where $\eta$ is as in \Cref{l-1} and $C$ depends only on $n$.

If $f\in BMO(B_1)$, then $u\in W^{2, BMO}(B_{1/2})$ and
\begin{equation}\label{e5}
\|u\|_{W^{2,BMO}(B_{1/2})}\leq C\left(\|u\|_{L^2(B_1)}+\|f\|_{BMO, B_{1}}\right).
\end{equation}
\end{theorem}
\proof Without loss of generality, we may assume that $\|u\|^*_{L^2(B_1)}\leq 1$, $\|f\|^*_{L^2(B_1)}\leq \delta$ and $|f|_{*,0}\leq \delta$, where $\delta$ is as in \Cref{l-1}. Otherwise, we may consider
\begin{equation*}
  v=\frac{\delta u}{\|u\|^*_{L^2(B_1)}+\|f\|^*_{L^2(B_1)}+|f|_{*,0}}.
\end{equation*}

For any $0<r<\eta^2$, there exists $m\geq 2$ such that $\eta^{m+1}\leq r<\eta^m$. By \Cref{t-1}, there exists a polynomial $P_{m-1}$ such that
\begin{equation*}
  \|u-P_{m-1}\|^*_{L^2(B_{\eta^{m-1}})}\leq \eta^{2(m-1)}.
\end{equation*}
Let $y=x/\eta^{m-1}$ and $v(y)=(u-P_{m-1})/\eta^{2(m-1)}$. Then
\begin{equation*}
  \Delta v=\tilde{f} ~~~~\mbox{in}~B_1,
\end{equation*}
where $\tilde{f}(y)=f(x)-f_{B_{\eta^{m-1}}}$. By the standard $W^{2,2}$ regularity for $v$,
\begin{equation*}
  \|D^2u-D^2P_{m-1}\|^*_{L^2(B_{\eta^{m}})}=\|D^2v\|^*_{L^2(B_{\eta})}
  \leq C\left(\|v\|^*_{L^2(B_{1})}+\|\tilde{f}\|^*_{L^2(B_{1})}\right)\leq C.
\end{equation*}
Hence,
\begin{equation}\label{e2.1}
 \|D^2u-(D^2u)_{B_r}\|^*_{L^2(B_{r})}\leq\|D^2u-D^2P_{m-1}\|^*_{L^2(B_{r})}
 \leq\frac{1}{\eta}\|D^2u-D^2P_{m-1}\|^*_{L^2(B_{\eta^{m}})} \leq C.
\end{equation}
That is, $u\in W^{2, BMO}(0)$.

If $f\in BMO(B_1)$ with radius $r_0$, $f\in BMO(x_0)$ for any $x_0\in B_{1/2}$ with the same radius $r_0$. Then by similar arguments to the above, we have $u\in BMO(x_0)$ for any $x_0\in B_{1/2}$ with the same radius $\eta^2r_0/2$. Hence, $u\in W^{2, BMO}(B_{1/2})$ and the estimate \cref{e5} holds.\qed~\\

\begin{remark}\label{r-1}
The techniques in \Cref{l-1} and \Cref{t-1} are adopted (with minor modifications) from \cite{lian2020pointwise} (see Section 11 there). The observation that \cref{e4} implies the BMO regularity is from \cite{MR3357696} (see Remark 6.3 there).
\end{remark}

The $W^{2,2}$ estimate for \cref{e1.1} is easy to establish (see \cite[Theorem 8.8, Theorem 9.9.]{MR1814364}). Since we have derived $W^{2, BMO}$ estimate, by the interpolation between $L^2$ and $BMO$ (see \cite[Appendix 4]{MR1616087} and \cite[Chapter 6.3.3]{MR3099262}), the $W^{2,p}$ estimate follows for any $2<p<\infty$. Furthermore, we can also obtain the $W^{2,p}$ estimate for any $1<p<2$ due to the duality (see \cite[Theorem 9.9.]{MR1814364}). In conclusion, we have the following Calder\'{o}n-Zygmund estimate:
\begin{corollary}\label{co-1}
Let $u\in W^{1,1}_{loc}(B_1)$ be a weak solution of
\begin{equation*}
\Delta u=f ~~~~\mbox{in}~B_1.
\end{equation*}
Suppose that $f\in L^{p}(B_1)$ with $1<p<\infty$. Then $u\in W^{2,p}_{loc}(B_1)$ and
\begin{equation}\label{e5}
\|u\|_{W^{2,p}(B_{1/2})}\leq C\left(\|u\|_{L^{1}(B_1)}+\|f\|_{L^{p}(B_1)}\right),
\end{equation}
where $C$ depends only on $n$ and $p$.
\end{corollary}

Based on the notion of the sharp maximal function $f^{\sharp}(x):=|f|_{\ast,x}$, we have another proof of the $W^{2,p}$ estimate without interpolation. Recall a property of the sharp maximal function due to Fefferman and Stein (see \cite[Theorem 6.30]{MR3099262}):
\begin{lemma}\label{le2.1}
Suppose that $g\in L^1(B_1)$ and $g^{\sharp}\in L^p(B_1)$ for some $1<p\leq \infty$. Then $g\in L^p(B_1)$ and
\begin{equation*}
\|g\|_{L^p(B_1)}\leq C\left(\|g^{\sharp}\|_{L^p(B_1)}+\|g\|_{L^1(B_1)}\right),
\end{equation*}
where $C$ depends only on $n$ and $p$.
\end{lemma}

In fact, we infer from the proof of \Cref{t-2} (see \cref{e2.1}) that
\begin{equation*}
(D^2u)^{\sharp}(x)\leq C\left(\|u\|_{L^2(B_1)}+\|f\|_{L^2(B_1)}+f^{\sharp}(x)\right),
~\forall ~x\in B_{1/2}.
\end{equation*}
Then for $2\leq p<\infty$,
\small
\begin{equation*}
  \begin{aligned}
\|D^2u\|_{L^p(B_{1/2})}
&\leq C\left(\|(D^2u)^{\sharp}\|_{L^p(B_{1/2})}+\|D^2u\|_{L^2(B_{1/2})}\right)\\
&\leq C\left(\|u\|_{L^2(B_1)}+\|f\|_{L^2(B_1)}+\|f^{\sharp}\|_{L^p(B_{1/2})}\right)\\
&\leq C\left(\|u\|_{L^2(B_1)}+\|f\|_{L^p(B_1)}\right),
  \end{aligned}
\end{equation*}
where $C$ depends only on $n$ and $p$.

\emph{Acknowledgements}. The authors would like to thank Professor Congming Li for a useful discussion.

\bibliographystyle{amsplain}
\bibliography{112}

\end{document}